\documentclass[11pt]{amsart}
\usepackage{amssymb}
\usepackage[
            pdftex,
           colorlinks=true,
            urlcolor=blue,       % \href{...}{...} external (URL)
            filecolor=blue,     % \href{...} local file
            linkcolor=blue,       % \ref{...} and \pageref{...}
            citecolor=blue,         % \cite{...}
            pdftitle= {Effective reconstruction of generic genus 4 curves from
theta hyperplanes},
            pdfauthor={David Lehavi},
            pdfsubject={},
            pdfkeywords={}
            pagebackref,
            pdfpagemode=None,
            bookmarksopen=true]
            {hyperref}
\usepackage{hyperref}
\usepackage{color}

\newcommand{\BP}{{\mathbb{P}}}
\newcommand{\BC}{{\mathbb{C}}}
\newcommand{\BZ}{{\mathbb{Z}}}

\newcommand{\gb}{\beta}

\newcommand{\gs}{\sigma}
\newcommand{\gS}{\Sigma}
\newcommand{\gom}{\omega}
\newcommand{\ga}{\alpha}
\newcommand{\gth}{\theta}
\newcommand{\gT}{\Theta}
\newcommand{\gt}{\tau}
\newcommand{\cO}{{\mathcal{O}}}
\newcommand{\cV}{{\mathcal{V}}}

\newcommand{\ti}[1]{\tilde{#1}}
\newcommand{\ol}[1]{\overline{#1}}

\newcommand{\Pic}{\mathrm{Pic}}

\newcommand{\Div}{\mathrm{Div}}
\newcommand{\sm}{\smallsetminus}

\theoremstyle{plain}
\newtheorem{lma}{Lemma}%[section]
\newtheorem{thm}[lma]{Theorem}
\newtheorem*{Thm}{Theorem}
\newtheorem{prp}[lma]{Proposition}
\newtheorem{cor}[lma]{Corollary}
\theoremstyle{definition}
\newtheorem{prd}[lma]{Proposition-Definition}
\newtheorem{rmr}[lma]{Remark}

\begin{document}
\title[Reconstruction of genus 4 curves from theta hyperplanes]{Effective reconstruction of generic genus 4 curves from their theta hyperplanes}
\begin{abstract}
Effective reconstruction formulas of a curve from its theta hyperplanes
are known classically in genus $2$ (where the theta hyperplanes
are Weierstrass points), and $3$ (where, for a generic curve, the theta
hyperplanes are bitangents to a plane quartic). However, for higher genera, no
formula or algorithm are known.
In this paper we give an explicit (and simple) 
algorithm for computing a generic genus $4$
curve from it's theta hyperplanes.
\end{abstract}
\author{David Lehavi}
\email{dlehavi@gmail.com}
\date{\today}
\maketitle
%%%%%%%%%%%%%%%%%%%%%%%%%%%%%
%
\section{Introduction}\label{Sintro}
%
%%%%%%%%%%%%%%%%%%%%%%%%%%%
The quest for methods of reconstructing a curve from its theta hyperplanes
goes back to
the 19th and early 20th century geometers Aronhold and Coble: in the non hyperelliptic genus
3 case, theta hyperplanes are simply
bitangents, and both Aronhold and Coble provided formulas for reconstructing curves
from certain ordered subsets of the $28$ bitangents of the curve (see
\cite{A}, \cite{Co} chapter IV, and \cite{Dol} sections 6.1.2 and 6.2.2).

Recent years witnessed some revived interest in generalizations of this
problem from several directions: first relaxing the need for {\em ordered}
theta hyperplanes (see \cite{CS1}, \cite{L}), and then generalizations to
higher genus curves (see \cite{CS2}) and abelian varieties (see \cite{GS-M1}, \cite{GS-M2}). However for $g>3$, these results are
not effective. In this work we give an effective result for the generic
genus $4$ case. Our result assumes ordered theta characteristic; that being
said, by \cite{CS2} this requirement is redundant for generic curves.

Throughout this paper we consider a generic complex curve $C$ of
genus $4$; since $C$ is generic we assume that all its odd theta
characteristics
are $1$ dimensional - i.e. if $\gth$ is an odd theta characteristic of $C$
then $\dim H^0(\gth,C)=1$. 
Hence, for each
odd theta characteristic $\gth$ there exists a unique hyperplane $l_\gth$
-- called a {\em theta hyperplane} -- in the
dual canonical system of $C$ such that when $C$ is identified with it's
canonical image, the points in the intersection product 
$C\cdot l_\gth$ are all double, and the points in $\frac{1}{2}C\cdot l_\gth$
sum up to $\gth$. This hyperplane is the projectivization of the plane
$T_\gth\gT_C\subset T_\gth JC$ under the identification of $T_\gth JC= T_0 JC
= H^0(K_C)^*$.
Recall (see e.g. \cite{Dol} 5.4.2), that if $\ga$ is a non-trivial
$2$ torsion point on the Jacobian $JC$, then the 
{\em Steiner system} $\gS_{C,\ga}$ of the pair $(C,\ga)$ is defined to be the 
set
\[
\{\gth: 2\gth = K_C\text{ and } \dim H^0(\gth,C)=\dim H^0(\gth+\ga,C)\equiv 1 \mod 2\}.
\]
The number of theta characteristics in a Steiner
system of a genus $g$ curve is $2^{g-1}\cdot(2^{g-1}-1)$; i.e. in our case
a Steiner system is comprised of $8\cdot7=56$ odd theta characteristics,
out of the total of $2^{4-1}(2^4-1)=120$
odd theta characteristics of the Jacobian. 
For each pair $\gth,\gth+\ga$, and corresponding theta hyperplanes
$l_\gth, l_{\gth+\ga}$, we let $q_{\{\gth,\gth+\ga\}}\in |\cO_{|K_C|^*}(2)|$ be the
image of $\{l_\gth, l_{\gth+\ga}\}\in S^2|K_C|$ under the embedding:
$S^2|K_C|\to |\cO_{|K_C|^*}(2)|$. 

We set $V_C:=H^0(\cO_{|K_C|^*}(2))$, and denote 
the map $H^0(\cO_{|K_C+\ga|^*}(2))\to H^0(2K_C)$ by $i$, 
the projection $V_C \to H^0(2K_C)$ by $p$,
and the pre-image $p^{-1}i H^0(\cO_{|K_C+\ga|^*}(2))\subset V_C$ by $V_{C,\ga}$.
Finally, we denote the
projectivization of a linear space by $\BP$.
The first component of our reconstruction algorithm is given by the
following: 
\begin{thm}\label{thm:QC}
Let $C,\ga,V_C,$ and $V_{C,\ga}$ be as above, (specifically, recall that
$C$ is generic), and let $Q_C$ be the unique
quadric surface containing the canonical image of $C$, then
\begin{enumerate}
\item\label{it:dim6} for all $\ga\in JC[2]\sm\{0\}$
we have
\[
 \text{span}(\{q_{\{\gth,\gth+\ga\}}\}_{\gth\in\gS_{C,\ga}})= \BP V_{C,\ga};
\]
\item\label{it:intersection}
moreover,
$\cap_{\ga\in JC[2]\sm\{0\}} \BP V_{C,\ga} = \{[Q_C]\}$,
where $[Q_C]$
denotes the moduli point representing $Q_C$ in the space $\BP V_C$.
\end{enumerate}
\end{thm}
Note that the theorem may be rephrased in the following ``genus free'' terms:
{\em the intersection
$\cap_{\ga\in JC[2]\sm\{0\}}\text{span}(\{q_{\{\gth,\gth+\ga\}}\}_{\gth\in\gS_{C,\ga}})$ 
is the locus of quadrics over $|K_C|^*$ containing $C$}. The beauty in this
statement
arises from its relation to the Enriques-Babbage Theorem (see e.g. 
\cite{ACGH} chapter VI \S3), which states that a canonical curve is either 
trigonal, or is isomorphic to a plane quintic, or is {\em cut out by quadrics}.

The proof of the theorem eventually reduces to the analysis of 
bi-elliptic curves: the first part is proved by considering a (degenerated)
double cover of the nodal cubic; whereas the second part is proved by
considering a curve with a big automorphism group.

Turning our attention back to the genus $4$ case, we now aim to locate
an irreducible cubic surface containing the canonical image of the curve $C$.
To this end we will make heavy use of the following classical Theorem:
\begin{Thm}[Wirtinger - see \cite{W}, \cite{Co} chapter V, or \cite{Ca}]
Let $C$ be a generic genus $4$ complex curve,
and let $\ga$ be a non trivial $2$-torsion point on the Jacobian of $C$.
Then the image of $C$ in the Prym canonical system $|K_C+\ga|^*$ is a sextic
with $6$ nodes, which are the intersection points of four lines in general
position: $l_1, l_2,l_3,l_4$ (not to be confused with the hyperplanes
$l_\gth$; we will never use the two notations in the same context). Moreover,
let $S$ be the blow-up of $|K_C+\ga|^*$ at the six intersection
points $l_i\cap l_j$, let $H$ be the pullback to $S$ of a generic hyperplane
in $|K_C+\ga|^*$ and let $E$ be the exceptional divisor in $S$,
then the complete linear system $|\cO_S(3H - E)|$ is canonically isomorphic
to $|K_C|$. Finally, the image of $S$ in the dual of this linear system
-- which we will denote below by $W_{C,\ga}$ -- is a Cayley cubic
(the unique space cubic with four nodes), where the four nodes are the blow
downs
of the strict transforms in $S$ of the four lines $l_1,\ldots, l_4$,
and where by its definition $W_{C,\ga}$ contains the canonical image of $C$.
\end{Thm}
Armed with this theorem we can state the following:
\begin{thm}\label{thm:WCa}
The four hyperplanes through each of the four
triplets of nodes of $W_{C,\ga}$ are --
set theoretically -- the four intersection points in $\BP V_C$ 
of the six dimensional projective subspace
\[
\BP\left(\left(V_C^*/((V_C/V_{C,\ga})^*\wedge(V_C/V_{C,\ga})^*)\right)^*\right),
\]
and the 2nd Veronese image of $|K_C|$ in  $\BP V_C$.
Moreover, each of these four points has multiplicity $2$ in the intersection.
\end{thm}
Given Theorems \ref{thm:QC} and \ref{thm:WCa},
the canonical image of the curve $C$ is readily reconstructed as the
intersection of $Q_C$
and $W_{C,\ga}$ -- which is the only cubic with nodes at the 
four intersection points of triplets of the four hyperplanes
determined in Theorem \ref{thm:WCa}.
%%%%%%%%%%%%%%%%
\subsection*{Acknowledgment}\ \\
Sam Grushevsky and an anonymous reviewer read early versions of this work
and gave the author extremely valuable feedback.
%%%%%%%%%%%%%%%%%%%%%%%%%%%%%
%%%%%%%%%%%%%%%%%%%%%%
%
\section{Proof of theorem \ref{thm:QC}}\label{sec:proofQC}
%
%%%%%%%%%%%%%%%%%%%%%%%%
We start by analyzing special Steiner systems of bi-elliptic curves.
Below we denote the ramification
(resp. branch) locus of a map $-$ by $R_{-}$
(resp. $B_{-}$).
\begin{thm}[Coble, \cite{Co} chapter IV]\label{thm:Coble}
Let $\pi:C_{be}\to E$ be a bi-elliptic double cover
where the genus of $C_{be}$ is $4$. We first note that $C_{be}$ is not
hyperelliptic.
Let $\ol{\pi}$ be the involution induced by $\pi$ on the dual
canonical system $|K_{C_{be}}|^*$ (which is a projective linear involution).
Then the fixed spaces of $\ol{\pi}$ are
a hyperplane $H_\pi$, and a point which we call the {\em focal point}.
Identifying $C_{be}$ with its canonical image, the intersection
$H_\pi\cdot C_{be}$ is the ramification divisor $R_\pi$,
which is comprised of six distinct points;
the image of $C$ in $H_\pi$ under
the projection from the focal point is $E$ embedded as a plane cubic;
$H_\pi$ is naturally
identified with a $g^2_3$ on $E$, denoted by $|L|^*$ which 
satisfied $2L\sim B_\pi$.

Conversely, the curve $C_{be}$ 
can be reconstructed from such data of $E$, $B_\pi$ (comprised of
six distinct points), and $L$ 
in the following way:
Construct $\BP^3$ as a cone over $|L|^*$, and in it reconstructs $C_{be}$
as the intersection of two surfaces:
\begin{itemize}
\item the cone $\ti{E}$ over the image of $E$ in the linear system $|L|^*$,
and,
\item a quadric surface $\ti{Q}$ in $\BP^3$ ramified over $|L|^*$ at the
unique conic satisfying intersecting $E$ at $B_\pi$.
\end{itemize}
\end{thm}
\begin{proof}
We first prove, arguing by contradiction, that
$C_{be}$ is not hyperelliptic. Suppose it is, and
denote the
hyperelliptic system on it by $|H_{C_{be}}|$. Since the hyperelliptic
involution commutes with all other involutions, we would have two double covers
$f:|H_{C_{be}}|^*\to\BP^1$ and $g:E\to\BP^1$ such that
$C_{be}=|H_{C_{be}}|^*\times_{\BP^1}E$. However, we then 
have the following inequality in $\mathrm{Div} E$:
\[
  B_\pi \leq f^* B_g ,
\]
which is impossible since $\deg B_g = 2\Rightarrow\deg f^*B_g=4$, while
$\deg B_\pi = 6$.

Note that since $\pi$ is a double cover, $B_\pi$ is a sum of six distinct
points.
We proceed to analyze the action of $\ol{\pi}$: As
$K_E$ is trivial, the image of $R_\pi$ in $|K_{C_{be}}|^*$ is cut out by some
hyperplane  $H_\pi$. Moreover, since $C_{be}$ does not admit a $g^1_2$, at most 
$3$ of the $6$ points of $R_\pi$ are collinear. Since $6$ points on $H_\pi$,
no $4$ of which are collinear, are fixed by $\ol{\pi}$,
the entire plane $H_\pi$ is also fixed by $\ol{\pi}$.

Recall that a projective linear involution which admits a fixed
hyperplane also admit a fixed point out of this hyperplane.
We call this point the {\em focal point} and project
the dual canonical system $|K_{C_{be}}|^*$ from this point to $H_\pi$.
As $\ol{\pi}$ is linear, all the lines through the focal
point (and some point on $H_\pi$) are $\ol{\pi}$ invariant. Hence, the
intersections of the canonical image of $C_{be}$ with these lines are
exactly the fibers of $\pi$, and the degree $2$ map $C_{be}\to H_\pi$ 
factors through a map $E\to H_\pi$; moreover, as $C_{be}$ is a degree $6$
embedding, the induced map $E\to H_\pi$ is a degree $3$ embedding.
Setting $L\in \Div^3 E$ to
be the (pullback to $E$ under an the embedding of the)
intersection of the image of $E$ in $H_\pi$ and some line in $H_\pi$,
we may identify $H_\pi$ with the dual complete linear system $|L|^*$.
Moreover, since the images of the $B_\pi$ in $|L|^*$
is the intersection of the image of $E$ in $H_\pi$, and the quadric surface
$Q_{C_{be}}$, they sit on a conic in $H_\pi$; thus, we have 
$2L\sim B_\pi$ in $\Pic^6 E$.
Expressing $|K_{C_{be}}|^*$
as a cone over $|L|^*$, we see that $C_{be}$ is the intersection of
 the cone over the image of $E$ in $|L|^*$ through the focal point, and the
quadric surface ramified over
$|L|^*$ at the unique conic passing through the images of $B_\pi$
there.
\end{proof}
Henceforth, we will identify $E$
with it's image in $|L|^*$, and $C_{be}$ with it's canonical image.
We now turn to the identification of some of the theta hyperplanes of $C_{be}$:
\begin{prd}\label{prd:24}
Assuming $C_{be}$ does not admit a theta null,
there are exactly $24$ theta hyperplanes of $C_{be}$ invariant under the
bi-elliptic involution. They are given as follows: Denote the $6$
distinct points of $B_\pi$ by $b_1,\ldots, b_6$; for each $b_i$
let $x_{i 1}, x_{i2}, x_{i3},x_{i4}$ be
the four points in $E$ satisfying $2x_{i j}+b_i \sim L$. 
For each $i,j$ we denote by $l_{i j}$ the line which satisfies
$l_{i j}\cdot E = 2x_{i j}+b_i$.
Let $H_{i j}$ be the pullback of $l_{i j}$ to $|K_{C_{be}}|^*$ under the identification
of $|K_{C_{be}}|^*$ as a cone over $|L|^*$, then
\[
  H_{i j}\cdot C_{be} = \pi^{-1}(2x_{i j}+b_i).
\]
\end{prd}
This proposition is an explicit form of proposition 2 in \cite{B}, where the genera of the curves
involved are $4$ and $1$.
\begin{proof}
Identifying $C_{be}$ with it's canonical image we have
$H_{i j}\cdot C_{be} = 2\pi^{-1}x_{i j}+\pi^{-1}b_i$; hence $H_{i j}$
represents an effective theta characteristic. Since we assume that $C_{be}$
does not admit a theta null, $H_{i j}$ is a theta hyperplane. 

Conversely, assume that $H$ is a theta hyperplane invariant under the bi-elliptic involution,
and let $l$ be the projection of $H$ to the linear system $|L|^*$, then the
following properties hold:
\begin{itemize}
\item If some $b_i$ satisfies $E\cdot l>b_i$, then $H\cdot C_{be} > \pi^{-1}(b_i)$,
which is the ramification point lying over $b_i$ -- with multiplicity $2$.
\item If some point $y\neq b_1\ldots,b_6$ satisfies $l\cdot E>a y$ for some positive $a$,
then both points in $\pi^{-1}(y)$ are in the intersection product $H\cdot C_{be}$,
 each with intersection multiplicity $a$.
\end{itemize}
Thus, if the intersection product $l\cdot E$ contains $B$ branching points and $n$ other points with positive
intersection multiplicities $a_k$, then
\[
  2(B+a_1+\cdots+a_n) = \# C_{be}\cdot H = 6,\quad\text{where all the $a_i$s are even}.
\]
The case where $n=0$ is the case
where $Q_{C_{be}}\cap H_\pi$ contains the line $l$,
which implies that this intersection is the
union of two lines, which implies that $Q_{C_{be}}$ is singular,
which implies that $C_{be}$ has a theta null.
Whence,
we have only one possible solution: $B=n=1,\quad a_1=2$.
\end{proof}
\begin{prp}\label{dsc:special}
Let $\gb$ be a non trivial $2$ torsion point in $J E$, then for each
$i=1,\dots,6$ and $j=1,2,3,4$ there is some $j'$ such that
$x_{i j}-x_{i j'}=\gb$.
This is a pairing on the $x_{i j}$s, which induces a pairing on the $H_{i j}$s.
Moreover, denoting by $\pi^*$ the map
$\Pic(E)\to\Pic(C_{be})$ induced by $\pi$, we have
$\pi^*x_{i j}- \pi^*x_{i j'}=\pi^*\gb\in JC_{be}[2]\sm\{0\}$.
\end{prp}
\begin{proof}
For each $i,j$ the shifts of $x_{i,j}$ by the four points of $J E[2]$
give the four $x_{i,j^\dagger}$, where $j^\dagger=1,2,3,4$.
Hence, $\gb$ induces a natural partition to pairs $x_{i j},x_{i j'}$, so that
$x_{i j}-x_{i j'}=\gb$ for all $i,j$. The second part follows since
$\pi^*$ is an embedding.
\end{proof}
In the proof of Corollary \ref{lma:Edim5} below, as well as in the proof of
Theorem \ref{thm:QC} we will apply a restricted form the following: 
\begin{prp}[Degeneracy loci of maps between vector bundles]\label{prp:degen}
The degeneracy loci of a map between two vector bundles over a base scheme
are closed subschemes of the base. 
\end{prp}
This proposition, as well as the proof, are classical. The proposition follows
from a  choice of a local basis to the bundles -- which 
is possible since the statement of the proposition is local, from
induction on the degeneracy rank, and
from the fact that the determinant is trivial on a
closed sub scheme.

As we just indicated, we will not apply the full strength of Proposition
\ref{prp:degen}, but rather the following corollary:
\begin{cor}\label{cor:span_dim}
Let $\cV/X$ be a vector bundle over a base $X$, and let
$\cV_1,\dots \cV_n$ be sub-bundles of $\cV$. 
Then the function $\dim \langle \cV_1|_x,\ldots, \cV_n|_x\rangle$ 
is lower semi-continuous on $X$,
and the function $\dim (\cap_{i=1}^n\cV_i|_x)$ is upper semi-continuous on $X$.
\end{cor}
\begin{proof}
By induction it suffices to prove this claim for two sub-bundles.
By Proposition \ref{prp:degen} the degeneracy locus
of the map $\cV_1\oplus\cV_2\to\cV$ is a closed subscheme of the base. Since the
degeneracy loci of both $\cV_1\to\cV$ and $\cV_2\to\cV$ are empty by definition
(as they are sub-bundles), we are done.
\end{proof}
\begin{lma}\label{lma:Edim5}
Let $E,L, \gb, b_i,\dots b_6$, and $x_{i j}$ be as above and generic, then 
the twelve reducible conics $l_{i j}\cup l_{i j'}$ in $|L|^*$
span the $5$ dimensional space $|\cO_{|L|^*}(2)|$.
\end{lma}
Note that Lemma \ref{lma:Edim5} has nothing to do with $C_{be}$; it is a
statement about plane cubics per se.
\begin{proof}
We apply the lower
semi-continuity of the dimension of span from Corollary \ref{cor:span_dim}.
We consider the degenerated case where $E$ is
a nodal cubic, and $\gb$ is degenerated to the trivial 2-torsion point.
As a model for the nodal cubic we will use the plane cubic $E_0$ given by the
nulls of $x^3+y^3-x y z$. The isomorphism between $\BC^*$ and $E_0$ is given by
\[
  \phi:t\mapsto (-t:t^2:1-t^3).
\]
The generic ``combinatorial'' scenario,
where we have $24$ points $x_{i j}$ coming in quartets indexed by the $j$
coordinates, where the differences
$x_{i j_1}-x_{i j_2}$ are the four 2-torsion points and where $x_{i j}-x_{i j'}=\gb$, 
degenerates in the $E_0$ case to the following scenario: We
have $12$ ``doubled'' points $x_{i j}$ which come in pairs indexed by the $j$
coordinate, where $x_{i1},x_{i2}=\phi(\pm t_i)$ for some $t_i$, and where
each $x_{i j}$ is {\em paired with itself}. As in 
the generic case the six points $b_i$ have to sit on a conic; however, as
for this degenerated case we have $b_i=\phi(t_i^2)$,
this constrain now translates to the easier constrain: $\prod_{i=1}^6t_i^2=1$.
Finally, instead of $24$ $l_{i j}$s we now have $12$,
where each one is ``trivially paired'' with itself. These $12$ $l_{i j}$s are
given by
\[
  \begin{aligned}
  T_{E_0}(\phi(t)) =& (3x^2-y z:3y^2-x z:-x y)|_{x=-t, y=t^2,z=1-t^3} \\
  =& (3t^2-t^2(1-t^3):3t^4+t(1-t^3):-t^3)\sim(2t+t^4:1+2t^3:-t^2),
  \end{aligned}
\]
for $t=\pm t_1,\dots, \pm t_6$.
We will show that if $t_1^2,\ldots,t_5^2$ are all distinct, then the span
of the ``doubled'' $T_{E_0}(\pm t_i)$ for $i=1,\ldots, 5$ is the entire
space $|\cO_{|\BP^2|^*}(2)|$.
 
Squaring the projective linear form $T_{E_0}(t)$ we get (using the
lexicographic order on degree 2 monomials):
\[
 (4t^2+4t^5+t^8:2t+5t^4+2t^7: -2t^3-t^6: 1+4t^3+4t^6:-t^2-2t^5:t^4).
\]
Let 
\[
  \begin{aligned}
 v_t :=  & (4t^2+4t^5+t^8,2t+5t^4+2t^7, -2t^3-t^6, 1+4t^3+4t^6,-t^2-2t^5,t^4) \\
=&(4t^2+t^8,5t^4,-t^6,1+t^6,-t^2,t^4)
+t(4t^4:2,-2t^2+2t^6,4t^2,-2t^5,0),
  \end{aligned}
\]
then our aim is to show that the evaluations of $v_t$
at $\pm t_1, \pm t_5$ span the $6$ dimensional affine space
$H^0(\cO_{|L|^*}(2))$.
Taking $t\neq 0,\infty$, the linear span of
$v_{t_i}, v_{-t_i}$ is equal to
the linear span of
$\frac{1}{2}(v_{t_i}+v_{-t_i}), \frac{1}{2}(v_{t_i}-v_{-t_i})$.
Denoting $s=t^2$, the vectors
$\frac{1}{2}(v_{t_i}+v_{-t_i}), \frac{1}{2}(v_{t_i}-v_{-t_i})$
are given by
\[
  (4s+s^4,5s^2,-s^3,1+s^3,-s,s^2),\quad
  t(2s^2,1,-s+s^3,2s,-s^2,0)
\]
respectively.
We now consider the matrix whose rows are the evaluations of
the vector $(4s+s^4,5s^2,-s^3,1+s^3,-s,s^2)$
at five distinct values of $s$. Performing
the column operations
\[
\mathrm{col}_1\mapsto \mathrm{col}_1+4\mathrm{col}_5,\quad
\mathrm{col}_4 \mapsto \mathrm{col}_4+\mathrm{col}_3,
\quad \mathrm{col}_2\mapsto \mathrm{col}_2-5\mathrm{col}_6
\] on 
this matrix we get the matrix whose rows are the
evaluations of the vector $(s^4,0,-s^3,1,-s,s^2)$. Except for the second
(trivial) column, and up to permutations and sign changes of the columns,
this matrix is a Vandermonde matrix; thus it is of degree $5$ for any
$5$ distinct values of $s$. Since the
rows of the last matrix span the kernel of the operator
$v\mapsto v(0,1,0,0,0,0)^t$, retracing the column operations we performed
on the original matrix, we see that the span of the rows is the kernel of the
operator
$v\mapsto v(0,1,0,0,0,-5)^t$.

Finally, as
\[
 t(2s^2,1,-s+s^3,2s,-s^2,0)(0,1,0,0,0,-5)^t = t\neq 0 \quad\text{generically},
\]
we see that for five distinct non zero values of $t^2$,
the {\em projective} space spanned by
the squares of the $10$ of the $l_{i j}$s is $5$ dimensional.
\end{proof}
\begin{cor}\label{cor:dim5}
Let $C_{be}, E, b_i, L, \gb, x_{i j}$ be generic as above.
Then the map
\[
  \begin{aligned}
  |\cO_{|L|^*}(2)| &\to \BP V_{C_{be}} \\
  [q] &\mapsto [\text{the cone over $q$ through the focal point}]
  \end{aligned}
\]
is an embedding into $\BP V_{C_{be},\pi^*\gb}$,
whose image is spanned by the twelve 
$q_{\{\gth,\gth+\pi^*\gb\}}$'s corresponding to the
partition, described in Proposition-Definition \ref{prd:24},
of the $24$ hyperplanes $H_{i j}$ to twelve pairs.
\end{cor}
\begin{proof}
The map is an embedding since $|K_{C_{be}}|^*$
is a cone over $|L|^*$. By Lemma \ref{lma:Edim5}, $|\cO_{|L|^*}(2)|$
is spanned by the twelve reducible conics $l_{i j}\cup l_{i j'}$; hence
by Proposition-Definition \ref{prd:24},
the image is spanned by the twelve $q_{\{\gth,\gth+\pi^*\gb\}}$s described above.
Finally, the twelve
$q_{\{\gth,\gth+\pi^*\gb\}}$'s lie in $\BP V_{C_{be},\pi^*\gb}$ by their definition;
hence, so does their span.
\end{proof}
\begin{rmr}[An alternative view of Corollary \ref{cor:dim5}]\label{rmr:LR}
As $Q_{C_{be}}$ contains the curve $C_{be}$, projecting
$\BP V_{C_{be},\pi^*\gb}$ from $[Q_{C_{be}}]$ 
we get $\BP V_{C_{be},\pi^*\gb}$. Composing this projection
on the embedding from Corollary \ref{cor:dim5},
we get an
isomorphism between $|\cO_{|L|^*}(2)|$ and
$\BP V_{C_{be},\pi^*\gb}$. In fact, more is true:
in \cite{LR} (Theorem 2.9
in the journal version, or 2.15 in the arxiv version) it is proved,
by a careful analysis of Coble's construction, that
there is a natural isomorphism between $|\cO_{|L+\gb|^*}|$ and
$|\cO_{|K_{C_{be}}+\pi^*\gb|^*}|$.
\end{rmr}
\begin{proof}[Proof of part \ref{it:dim6} of Theorem \ref{thm:QC}]
In Corollary \ref{cor:dim5} we showed
that for $C_{be}, E,\gb$ as in the corollary there are exactly
twelve quadrics $q_{\gth,\gth+\pi^*\gb}$ passing through the focal point, and that
these quadrics span a five dimensional subspace
of $\BP V_{C_{be},\pi^*\gb}$. As the Steiner
system $\gS_{C_{be},\pi^*\gb}$ admits $56-24=32$ theta hyperplanes which
{\em do not}
pass through the focal point of the involution, there are
$16$ quadrics $q_{\gth,\gth+\pi^*\gb}$ which do not sit on the $5$ dimensional
space spanned by the $12$ quadrics from \ref{prd:24}. Hence, 
$\dim\text{span}(\{q_{\{\gth,\gth+\pi^*\gb\}}\}_{\gth\in\gS_{C_{be},\pi^*\gb}})\geq 6$.
By Corollary \ref{cor:span_dim} we see that
$\dim \text{span}(\{q_{\{\gth,\gth+\ga\}}\}_{\gth\in\gS_{C,\ga}})\geq 6$
for generic curves, and not only
for bi-elliptic ones.
However, as
$\dim \text{span}(\{q_{\{\gth,\gth+\ga\}}\}_{\gth\in\gS_{C,\ga}})\subset\BP V_{C,\ga}$,
which is $6$ dimensional,
we see that for a generic curve $C$, $\dim \BP V_{C,\ga}=6$.
\end{proof}
To prove the second part of theorem \ref{thm:QC} we analyze the
mutual structure of several Steiner systems on one curve. As we
have already analyzed
(some of) the structure of a specific Steiner system on a bi-elliptic curve, we
will introduce and analyze below a special genus $4$ curve with many
bi-elliptic involutions; this would immediately give us many Steiner systems
on this curve, of the form already analyzed above.
\begin{prd}[Kuribayahi and Kuribayashi, see \cite{KK} Proposition 2.4(f)(1)]\label{prd:many_inv}
There is exactly one genus $4$ curve,
denoted here by
$C_{9\times 8}$, whose automorphism group is isomorphic to $(\BZ/3)^2\ltimes D_8$.
Endowing the dual canonical system of $C_{9\times 8}$ with the coordinates
$x_1,\ldots x_4$, and
presenting the automorphism group as automorphism of the canonical system,
the group is generated by the automorphisms:
\[
  \begin{aligned}
   \gt_1=(x_1:x_2:x_3:x_4)&\mapsto (\gom x_1: \gom^2 x_2:x_3:x_4),
\quad\text{where }\gom^3=1,\\
   \gt_2=(x_1:x_2:x_3:x_4)&\mapsto (-x_2: -x_1:x_3:x_4),\quad\text{and}\\
   \gt_3=(x_1:x_2:x_3:x_4)&\mapsto (x_3:x_4:x_1:x_2).
  \end{aligned}
\]
\end{prd}
\begin{proof}
See \cite{KK}. Verifying that the automorphism group contains
a group isomorphic to $(\BZ/3)^2\ltimes D_8$, which is all we need
for our purpose here,
is immediate using the explicitly model in Corollary \ref{cor:C98} below.
\end{proof}
\begin{cor}[Swinarski, see \cite{S}]\label{cor:C98}
The canonical image of the curve $C_{9\times 8}$, in the coordinates
of \ref{prd:many_inv}, is the intersection 
of the following cubic and quadric:
\[
  x_1^3-x_2^3+x_3^3-x_4^3,\qquad x_1x_2+x_3x_4.
\]
\end{cor}
\begin{proof}
Direct verification.
\end{proof}
\begin{prp}\label{prp:c98cover}
The quotient of $C_{9\times 8}$ under $\tau_2$ is elliptic.
Denoting the quotient by $\tau_2$ by $\pi_{\tau_2}$ and using our notations
from Theorem \ref{thm:Coble},
the focal point is $(1:1:0:0)$ and
 $H_{\pi_{\tau_2}}=\{(a:-a:b:c)|a,b,c\in \BC\}$. The image of the
elliptic quotient by $\tau_2$ in $H_{\pi_{\tau_2}}$ is given by the null set
of $(2a)^3+b^3-c^3-3(2a)b c$.
\end{prp}
\begin{proof}
The claim about the invariant sub-spaces follows immediately by direct
verification.
By Corollary \ref{cor:C98}, $C_{9\times 8}$ sits on the null set of the cubic
surface
\[
  x_1^3-x_2^3+x_3^3-x_4^3 - 3(x_1-x_2)(x_1x_2+x_3x_4) =
   (x_1-x_2)^3 +x_3^3-x_4^3 -3(x_1-x_2)x_3x_4.   
\]
The plane cubic curve is the intersection of this cone and the
invariant hyperplane; i.e. in the plane
coordinates given above it is the curve
\[
  (2a)^3+b^3-c^3-3(2a)b c.  
\]
\end{proof}
\begin{proof}[Proof of part \ref{it:intersection} of Theorem \ref{thm:QC}]
By definition, for any  $C,\ga$ we have 
\[
[Q_C]\in \cap_{\ga\in JC\sm\{0\}} \BP V_{C,\ga}.
\]
As 
$V_{C,\ga}$ are sub-bundles
of $V_C$ over the moduli of $C,\ga$,
Corollary \ref{cor:span_dim} tells us that in order to proof part
\ref{it:intersection} of Theorem \ref{thm:QC} we need only 
prove it for the curve $C_{9\times 8}$.

If $\gs$ is an involution of $C_{9\times 8}$, we will denote by $\ga_\gs$ one of the three $2$-torsion points
in $JC_{9\times 8}[2]\sm\{0\}$ invariant under the involution: these
are the pullback of the three non trivial 2-torsion points of the quotient by
$\gs$, which is an elliptic curve 
(``morally'', we don't care which one of the three we pick is that as we are
dealing with second symmetric products of linear systems).
By Corollary \ref{cor:dim5}, (and using the notation
$L$ as in the corollary),
the space $V_{C_{9\times 8},\ga_{\gt_2}}$ is spanned by
the pullback of $S^2H^0(L)$ to $S^2H^0(K_{C_{9\times 8}})=V_{C_{9\times 8}}$
under the
presentation of $|K_{C_{9\times 8}}|^*$ as a cone over $|L|^*$, and
a quadric form in $V_{C_{9\times 8},\ga_{\gt_2}}$ which is not trivial on the focal
point;
as the quadric $Q_{C_{9\times 8}}$ does not pass through the focal point,
and is in $V_{C_{9\times 8},\ga_{\gt_2}}$, we may choose this quadric form
to be $x_1 x_2 + x_3 x_4$.

In the 
notations used in \ref{prd:many_inv},
$L$ is spanned by $x_1-x_2, x_3, x_4$; hence:
\[
  \begin{aligned}
V_{C_{9\times 8},\ga_{\gt_2}} =&
  \langle x_1x_2+x_3x_4, (x_1-x_2)^2, (x_1-x_2)x_3, (x_1-x_2)x_4, 
   x_3^2,  x_3x_4, x_4^2\rangle \\
 = & \langle x_1^2+x_2^2, (x_1-x_2)x_3, (x_1-x_2)x_4,
   x_1x_2, x_3^2, x_3x_4, x_4^2\rangle.
  \end{aligned}
\]
By symmetry,
\[
  V_{C_{9\times 8},\ga_{\gt_2\gt_1}} =
\langle x_1^2+\gom^2x_2^2, (x_1-\gom x_2)x_3, (x_1-\gom x_2)x_4,
   x_1x_2, x_3^2, x_3x_4, x_4^2\rangle.
\]
In order to compute the intersection of these spaces,
we observe that $V_{C_{9\times 8}}$ may be
broken to the direct sum of the following spaces:
\[
   \langle x_1^2,x_2^2\rangle\oplus\langle x_1x_3, x_2x_3\rangle
    \oplus \langle x_1x_4,x_2x_4\rangle\oplus
  \langle x_1x_2, x_3^3,x_3x_4, x_4^2\rangle.
\tag{$\diamond$}
\]
Considering the spanning elements (indeed - bases) we chose for
the spaces $V_{C_{9\times 8},\ga_{\gt_2}}$
and $V_{C_{9\times 8},\ga_{\gt_2\gt_1}}$, it is clear that
both spaces are direct sums of their intersections with the components in
equation ($\diamond$), that
both intersections with the last component in equation
($\diamond$) form the entire last component,
and finally that the intersections of these spaces with each of the
other three components in equation ($\diamond$) are
one dimensional and different from each other.
Hence we have
\[ 
 V_{C_{9\times 8},\ga_{\gt_2}}\cap V_{C_{9\times 8},\ga_{\gt_2\gt_1}}
  =  \langle x_1x_2, x_3^2, x_3x_4, x_4^2\rangle
  =  \langle x_1x_2+x_3x_4, x_3^2, x_3x_4, x_4^2\rangle.
\]
To conclude the proof, note that the intersection of this space with it's
image under $\gt_3$ is generated just by $x_1x_2+x_3x_4$, and recall that
the space generated by $x_1x_2+x_3x_4$ is by definition in the
intersection of all the $V_{C_{9\times 8},\ga}$'s.
\end{proof}
%%%%%%%%%%%%%%%%%%%%%%%%%%%%%%%%%%%%%%%%%%%%%%%%%%%%%%%
%
\section{Proof of theorem \ref{thm:WCa}}\label{sec:proofWCa}
%
%%%%%%%%%%%%%%%%%%%%%%%%%%%%%%%%%%%%%%%%%%%%%
\begin{lma}\label{lma:WCa}
Let $X$ be the intersection of
the image of $|K_C|\times |K_C|$ in $\BP V_C$ with
$\BP V_{C,\ga}$, then
the two dimensional fibers of the projection on the first
coordinate of the pullback of $X$ to $|K_C|\times |K_C|$ lie
exactly over the moduli points which represent the four hyperplanes through
the four possible triplets of the blow-downs
of the strict transforms in $S$ of the four lines $l_1,\ldots,l_4$
(see Wirtinger's Theorem in the introduction regarding $S$ and the $l_i$s).
\end{lma}
\begin{proof}
Note that the isomorphism $|\cO_{S}(3H-E)|\cong |K_C|$ from Wirtinger's theorem
induces an isomorphism between quadric forms 
over $|K_C|^*$, and sextic forms over $|K_C+\ga|^*$ containing the six
intersection points $l_i\cap l_j$.
Moreover, if we let $L_i$ be a linear form on $|K_C+\ga|^*$
whose null set is $l_i$, then this map sends any quadric form on $|K_C|^*$
which
is zero on the four nodes of $W_{C,\ga}$ to a product of the four $L_i$'s and
a quadric in $|K_C+\ga|^*$.

We write the cubic $W_{C,\ga}$ as the null set:
\[
  W_{C,\ga} = Z(h_1h_2h_3+h_1h_2h_4 + h_1h_3h_4+h_2h_3h_4),
\]
where $h_1, h_2, h_3, h_4$ are linear forms such that 
the four points $\cap_{1\leq i \leq 4, i\neq j}Z(h_i)$ -- for
$j\in\{1,2,3,4\}$ -- are the four nodes of $W_{C,\ga}$.  We endow the space
$V_C$ with the coordinates $\{h_i h_j\}_{1\leq i \leq j\leq4}$.
Then the quadrics passing through the four nodes of $W_\ga$ are exactly 
the null
sets of quadric forms in the span of the six forms:
$\{h_i h_j\}_{1\leq i < j\leq4}$.

Recall that the curve $C$ is generic on the moduli of genus $4$
curves, which is $9$ dimensional. Hence, writing $Q_C$ as the null-set of
$q_C=\sum_{1\leq i<j\leq 4} d_{i j}h_i h_j+\sum_{i=1}^4 a_i h_i^2$
(not to be confused with the $q_{\gth,\gth+\ga}$ which already appeared before),
$q_C$ has ten free parameters; thus w.l.o.g. we may assume that all the
$a_i$s are non 0. As
$q_C$ is not in the span of $\{h_i h_j\}_{1\leq i<j\leq 4}$, the 
space $V_{C,\ga}$ is spanned inside
$V_C$ by $q_C$ and $\{h_i h_j\}_{1\leq i<j\leq 4}$.
We can now compute the intersection $\BP V_{C,\ga}$
with the image of the map:
\[
  s:|K_C|\times|K_C|\to \BP V_C.
\]
Let $\sum_{i=1}^4 b_i h_i$ and $\sum_{i=1}^4 c_i h_i$ be two canonical sections,
then their product modulo $\mathrm{span}(\{h_i h_j\}_{1\leq i<j\leq 4})$ is 
$\sum_{i=1}^4 b_i c_i h_i^2$.
Hence, in order for the projectivization of the product of these canonical
sections to lie in $\BP V_{C,\ga}$ we must have an
equality
\[
  (b_1c_1,b_2c_3,b_3c_3) = k(a_1,a_2,a_3,a_4)\quad\text{for some }k.
\]
If all the $b_i$s are non $0$, then the only non trivial
solution to this
equality is the one where $k\neq 0$, and $c_i = k a_i/b_i$. 
If at least one $b_i$ is $0$ then $k$ is necessarily also $0$, and in this
case the indices for which the $b_i$s are $0$ are the one complimentary (in
the set $\{1,\ldots,4\}$) to the ones for which the $c_i$s are $0$. 

Pulling the image of the intersection back to $|K_C|\times |K_C|$ (under 
the map $s$), 
we see that the projection on the first coordinate is a birational
isomorphism between the intersection and $|K_C|$, where the exceptional fibers are of two types:
\begin{itemize}
\item The fibers over each of the moduli points $[h_i]$ is a $\BP^2$.
\item The fiber over each point in the six lines $span([h_i],[h_j])$, is a $\BP^1$.
\end{itemize}
\end{proof}
\begin{proof}[Proof of theorem \ref{thm:WCa}]
By Lemma \ref{lma:WCa} the four hyperplanes through
triplets of the four nodes of $W_{C,\ga}$ are exactly the
projective points represented by $a\in H^0(K_C)$, such 
that for all non trivial
$v\in (V_C/V_{C,\ga})^*\subset V_C^*=S^2(H^0(K_C))$,
the null spaces of the operators
\[
  \begin{aligned}
  H^0(K_C)&\to \BC \\
  x &\mapsto v\cdot (a\otimes x),\qquad 
     \text{where $\cdot$ is tensor contraction}
  \end{aligned}
\]
are identical.
Equivalently, this means 
that any two vectors in $(V_C/V_{C,\ga})^*\cdot a$ are dependant
I.e. the solutions are the solutions of the equation:
\[
((V_C/V_{C,\ga})^* \cdot a)\wedge((V_C/V_{C,\ga})^*\cdot a) = 0,
\]
However, since contraction commutes with tensor products (and wedge products,
which may be viewed as tensor products followed by projections to a subspace), 
\[
 ((V_C/V_{C,\ga})^*\cdot a)\wedge((V_C/V_{C,\ga})^* \cdot a)=
 ((V_C/V_{C,\ga})^*\wedge(V_C/V_{C,\ga})^*) \cdot (a\otimes a).
\]
Hence the solution set is the intersection of
\[
\BP((V_C^*/((V_C/V_{C,\ga})^*\wedge(V_C/V_{C,\ga})^*))^*),
\]
and the image of $|K_C|^*$ under the 2nd Veronese map. Finally,
Since the degree of the 2nd Veronese of a $3$ dimensional space
is $2^3=8$, and since by Theorem C in \cite{Ca}, the Galois group
acting on the four nodes of $W_{C,\ga}$ is $S_4$; hence, each of the
intersection points has multiplicity $2$.
\end{proof}
%%%%%%%%%%%%%%%%%%%%%%%%%%%%%%%%%%%%%%%%%%%%%%%%%%%%%%%%%%%%%%%%%%%
%                                                                 %
% bibliography                                                    %
%                                                                 %
%%%%%%%%%%%%%%%%%%%%%%%%%%%%%%%%%%%%%%%%%%%%%%%%%%%%%%%%%%%%%%%%%%%

\end{document}